\documentclass[12pt,oneside,reqno]{amsart}
\usepackage{graphicx}
\usepackage{caption,subcaption}
\usepackage{dcolumn}
\usepackage{amsmath,amsthm,amssymb}
\usepackage{todonotes}
\usepackage{url}
\usepackage{placeins}

\theoremstyle{remark}
\newtheorem{rem}{Remark}

\begin{document}

\title[Generic self-similar blowup for equivariant wave maps]{Generic self-similar blowup for equivariant wave maps and Yang-Mills fields in higher dimensions}

\author{Pawe\l {} Biernat}
\address{Institute of Mathematics, Jagiellonian University, Krak\'ow, Poland}
\email{pawel.biernat@gmail.com}
\author{Piotr Bizo\'n}
\address{Institute of Physics, Jagiellonian
University, Krak\'ow, Poland\\ and
Max Planck Institute for Gravitational Physics (Albert Einstein Institute),
Golm, Germany}
\email{piotr.bizon@aei.mpg.de}

\date{\today}%
\begin{abstract}
 We consider equivariant wave maps from  the $(d+1)$--dimensional Minkowski spacetime into the $d$-sphere for $d\geq 4$. We find  a new explicit stable self-similar solution and give numerical evidence  that it plays the role of a  universal attractor for generic blowup. An analogous result is obtained for  the $SO(d)$ symmetric Yang-Mills field for $d\geq 6$.
\end{abstract}
\maketitle

\section{Introduction}
 This
paper is concerned with the Cauchy problem for the semilinear radial wave equation in $(d+1)$ dimensions
\begin{equation}\label{eq}
u_{tt} = u_{rr}+\frac{d-1}{r}\,u_r -\frac{f(u)}{r^2}\,,
\end{equation}
where $r\geq 0$ is the radial variable and  $u=u(t,r)$. For concretness, we shall focus on two specific and particularly interesting cases:
\begin{equation}\label{f}
f(u)=
\begin{cases}
(d-1) \sin(u) \cos(u)\,, & \text{(WM)} \\
d\, u (1-u) (2-u)\,, & \text{(YM)}
\end{cases}
\end{equation}
which correspond to equivariant wave maps from the $(d+1)$--dimensional Minkowski spacetime into the $d$-sphere and an $SO(d+2)$ symmetric Yang-Mills field in the $(d+3)$--dimensional Minkowski spacetime, respectively
 (see e.g. \cite{cst} for the derivation).

The basic question for equation \eqref{eq} is whether
solutions starting from  smooth initial data
 can become singular (``blow up'') in  finite time and, if so, how does the blowup occur. The key feature, which plays a major role in answering this question,  is the scale invariance of Eq.\eqref{eq} with respect to  the scaling tranformation  $u(t,r) \mapsto
u_{L}(t,r)=u(t/L,r/L)$, where $L$ is a positive constant. Under this scaling the
conserved
 energy
\begin{equation}\label{energy}
E(u)= \int\limits_0^{\infty} \left(u_t^2+ u_r^2+
  \frac{F(u)}{r^2} \right) \: r^{d-1} dr\,,
\end{equation}
where $F'(u)=f(u)$, transforms as
$E(u_{L})=L^{d-2} E(u)$, which means that Eq.\eqref{eq}
 is energy critical for $d=2$ and supercritical for $d\geq 3$.

 The critical case $d=2$ has been intensively studied for the past decades leading to good understanding of the mechanism of blowup. Most notably, Struwe proved that singularities necessarily have the form of self-shrinking harmonic maps (hence nonexistence of harmonic maps implies global regularity) and the rate of shrinking is faster than self-similar \cite{struwe}. For the nonlinearities \eqref{f} the generic rates of shrinking were derived formally in \cite{bos, os} and proved  rigorously in \cite{rr}.

    The understanding of the supercritical case $d\geq 3$ is much less satisfactory. All known examples of blowup involve self-similar solutions \cite{s,ts, cst, b1,b2}, however  it is not known if singularities \emph{must} be self-similar.  For $d=3$  and the nonlinearities \eqref{f}, the  self-similar blowup was proved (under a certain spectral condition) to be stable \cite{d1,d2} and demonstrated numerically to be generic \cite{bct,bt}. In higher dimensions, as far as we know, nothing was known in even dimensions while for odd $d\geq 5$ some self-similar solutions have been constructed \cite{cst,b1}  but their stability and role in dynamics have not been studied.

 In this paper we present new explicit self-similar solutions of Eq.\eqref{eq} in all dimensions $d\geq 4$ (section~2), analyze their linear stability (section~3), and finally give numerical evidence that they play the role of universal attractors in the generic blowup (section~4).
\section{Self-similar solutions}
By definition, self-similar solutions are invariant under the scaling, that is $u(t/L,r/L)=u(t,r)$, hence they have the form
\begin{equation}\label{css-ansatz}
  u(t,r)=\phi(y)\,,\qquad y=\frac{r}{T-t}\,,
\end{equation}
where a positive constant $T$, clearly allowed by the time translation symmetry, is introduced for later convenience.
Inserting this ansatz into Eq.\eqref{eq} one obtains the ordinary differential equation
\begin{equation}\label{css-ode}
(1-y^2) \phi'' +\left(\frac{d-1}{y}-2y\right) \phi'
-\frac{f(\phi)}{y^2}=0\,.
\end{equation}
We are interested in solutions that are smooth on the closed interval $0\leq
y\leq 1$, which corresponds to the interior of the past light
cone of the  point $(t=T,r=0)$. For such solutions
\begin{equation}\label{css-blow}
  \frac{\partial^{n}}{\partial r^n}\, \phi\left(\frac{r}{T-t}\right )\Big\vert_{r=0} = (T-t)^{-n} \, \frac{d^n\phi}{dy^n} (0)\,,
\end{equation}
hence the $n$-th derivative  at the origin (such that $\frac{d^n\phi}{dy^n} (0)\neq 0$) diverges as $t\nearrow T$. Thus, each self-similar solution $\phi(y)\in C^{\infty}[0,1]$ is an example of a singularity developing in finite time from smooth initial data.

We found the following explicit solutions of Eq.\eqref{css-ode}
\begin{equation}\label{phi0}
\phi_0(y)=
\begin{cases}
2 \arctan\left(\dfrac{y}{\sqrt{d-2}}\right)\,, & \text{(WM)} \\
\dfrac{a y^2}{y^2+b}\,, & \text{(YM)}
\end{cases}
\end{equation}
where
$$
a = \frac{2(d+4-\sqrt{3d(d-2)})}{3d-2\sqrt{3 d(d-2)}},\quad b = \frac{1}{3} \frac{\sqrt{3d(d-2)}\,(8-d)}{3d-2\sqrt{3 d(d-2)}}\,.
$$
For $d=3$ these solutions  have been known \cite{ts,b2}, however for $d\geq 4$ they appear to be new. In particular, the only smooth self-similar solutions obtained previously, for $d=5$ (WM) and $d\in \{5,7\}$ (YM), by  variational \cite{cst} and shooting \cite{b1} methods are different from $\phi_0$.  For both these methods it was essential that $\phi(1)=\phi_*$, where $\phi_*$ is the convexity radius of the target manifold, equal to $\pi/2$ for (WM) and 1 for (YM). In contrast, for our solutions $\phi_0(1)<\phi_*$ when $d\geq 4$.

  \begin{rem} A complete classification of smooth solutions of Eq.\eqref{css-ode} for $d\geq 4$ is not an easy problem; the main difficulty being due to an involved smoothness condition at $y=1$ (see section 2.2 in \cite{cst} for the analysis of this condition in odd dimensions). Fortunately, this hard ODE problem need not concern us here because, as  will be shown in section~4, smooth solutions other than $\phi_0$ do not seem to participate in the generic blowup, presumably due to their instabilities.
 \end{rem}

\begin{rem} Clearly, solutions similar to $\phi_0$ exist for a variety of nonlinearities other than \eqref{f} (albeit not in an explicit form, of course). In particular,  by modifying \eqref{f} for $\phi>\phi_0(1)$ one can construct analogous smooth solutions for convex targets using the shooting method of \cite{cst}.
\end{rem}

\begin{rem}
It is instructive to compare the problem of blowup for Eq.\eqref{eq} with the analogous problem for the heat flow
\begin{equation}\label{heat}
u_{t} = u_{rr}+\frac{d-1}{r}\,u_r -\frac{f(u)}{r^2}\,,
\end{equation}
with the same nonlinearities \eqref{f}. This equation has self-similar blowup solutions  for $3\leq d\leq 6$ (WM) \cite{fan} and $3\leq d\leq 7$ (YM) \cite{wein} but no such solutions exist in higher dimensions \cite{bw}. In accordance with that, it was found  that  the blowup for Eq.\eqref{heat} changes character from self-similar \cite{bb} to non-self-similar \cite{bier} in sufficiently high dimensions.
\end{rem}

\section{Linear stability analysis}
 As the first step towards understanding the role of the self-similar solution $\phi_0(\frac{r}{T-t})$ in dynamics we need to analyze its linear stability. To this end
 it is convenient to define a new time coordinate
$s=-\ln(T-t)$ and rewrite Eq.\eqref{eq} in terms of
$U(s,y)=u(t,r)$
\begin{equation} \label{eq-ys}
U_{ss} + U_{s} + 2 y\: U_{sy}
=(1-y^2) U_{yy} +\left(\frac{d-1}{y}-2y\right) U_{y}
-\frac{f(U)}{y^2}\,.
\end{equation}
In these variables the problem of finite time blowup in converted
into the problem of asymptotic convergence for $s \rightarrow
\infty$ towards the stationary solution $\phi_0(y)$.
 Following the standard procedure we seek
solutions of Eq.\eqref{eq-ys} in the form
$U(s,y)=\phi_0(y)+ e^{\lambda s} v(y)$. Dropping nonlinear terms
 we get the quadratic eigenvalue equation
\begin{equation}\label{eigen}
(1-y^2) v''+\left(\frac{d-1}{y}-2(\lambda+1)y\right) v' -\lambda(\lambda+1) v - V(y) v=0,
\end{equation}
where
\begin{equation*}\label{poten}
    V(y) := \frac{f'(\phi_0(y))}{y^2} =
    \begin{cases}
    \dfrac{d-1}{y^2}\,\dfrac{y^4+(12-6d)y^2+(d-2)^2}{(y^2+d-2)^2}\,, & \text{(WM)}\\
    \dfrac{d}{y^2} \, \dfrac{(3a^2-6a)y^4+(2-6ab)y^2+2b^2}{(y^2+b)^2} & \text{(YM)}\,.
    \end{cases}
\end{equation*}
We consider the eigenvalue problem \eqref{eigen} on the interval $0\leq
y\leq 1$. By assumption,  the solution $U(s,y)$ is smooth for $s<\infty$, hence we
 we demand that  $v\in \mathcal{C}^{\infty}[0,1]$. This condition leads to the quantization of the eigenvalues.

 In order to determine the analyticity properties of solutions of Eq.\eqref{eigen}
it is convenient to introduce new variables $x$ and $w(x)$ defined by
\begin{equation}\label{subs}
    x=\frac{
    c y^2}{y^2+c-1},\qquad v(y)= x^{\mu}\,(c-x)^{\frac{\lambda}{2}}\, w(x),
\end{equation}
where $\mu=1/2$,  $c=d-1$ for (WM) and  $\mu=1$, $c=b+1$ for (YM).
This change of variables brings Eq.\eqref{eigen} into the canonical form of the Heun equation \cite{nist}
\begin{equation}\label{heun-eq}
  w''+\left(\frac{\gamma}{x}+\frac{\delta}{x-1}+\frac{\epsilon}{x-c}\right)\, w'+\frac{\alpha\beta x -q}{x(x-1)(x-c)}\,w=0\,,
\end{equation}
where in the case of (WM)
\begin{eqnarray*}
\gamma &=&1+d/2, \quad \delta=(3-d)/2+\lambda,\quad \epsilon=d/2, \quad \alpha=(\lambda+5)/2, \\
 \beta &=&(\lambda-1)/2,\quad q=(\lambda-1)(d\lambda+5d-2\lambda-6)/4\,,
 \end{eqnarray*}
 and in the case of (YM)
 \begin{eqnarray*}
\gamma &=&2+d/2, \quad \delta=(3-d)/2+\lambda,\quad \epsilon=1/2,\quad \alpha=3+\lambda-\beta, \\
 \beta &=& \frac{1}{2b+2}\,\left(b\lambda+(3 a^2 b d+3 a^2 d+b^2+2b+1)^{\frac{1}{2}}+3b+\lambda+3\right),\\
 q &=& b \lambda^2/4-3a d/2+5 b \lambda/4+ d \lambda/4+3b/2+d/2+\lambda+2\,.
 \end{eqnarray*}
 Equation \eqref{heun-eq} has four regular singularities at $0,1,c,\infty$. The characteristic exponents at $0$ and $1$ are $\{0,1-\gamma\}$ and $\{0,1-\delta\}$, respectively.
The regular solutions at $x=0$ and $x=1$ are given by the local Heun functions (using the Maple notation)
\begin{eqnarray}\label{heun0}
  w_0(x)&=&\mathrm{HeunG}\left(c,q,\alpha,\beta,\gamma,\delta; x\right)\,,\\
  w_1(x)&=&\mathrm{HeunG}\left(1-c,-q+\alpha\beta,\alpha,\beta,\delta,\gamma; 1-x\right)\,.
\end{eqnarray}
The eigenvalues are given by zeros of the Wronskian
\begin{equation}\label{wron}
W[w_1,w_2](\lambda):=w_1'(x) w_2(x)-w_2'(x) w_1(x)\,.
\end{equation}
Unfortunately,  the connection problem for the Heun equation is unresolved and  therefore (in contrast to the hypergeometric equation) no explicit formula for the Wronskian is available. Fortunately, the Heun functions  are built into Maple and the zeros of the Wronskian can be computed numerically with great precision\footnote{An alternative way is to use the method of continued fractions \cite{b3}. We have verified that these two methods give identical results.}. The results of this computation are displayed in Table~1 (WM) and Table~2 (YM).

We remark that the  eigenvalue $\lambda_0=1$ is due to the freedom of changing the blowup time $T$ and does not correspond to any instability. The fact that the corresponding eigenfunction $v_0(y)=y \phi_0'(y)$ has no zeros can be used to exclude rigorously eigenvalues with $\Re(\lambda)>d-2$ (see \cite{b3,dsa} for the proof in $d=3$), however it seems hard to prove the nonexistence of eigenvalues with $0<\Re(\lambda)<d-2$. Nonetheless, we feel confident that if such an eigenvalue existed, we would have  found it numerically. Consequently, we assume that all the eigenvalues but $\lambda_0$ have negative real part\footnote{It is somewhat surprising that all eigenvalues are real.}. From this assumption it follows that the solution $\phi_0$ is linearly asymptotically stable \cite{dsa}.

\begin{table}[!h]
  \centering
   \scalebox{0.9}{
\begin{tabular}{|cc|ccccc|}\hline
\noalign{\smallskip}
&  & $n=0$ & $n=1$ & $n=2$ & $n=3$ & $n=4$ \\
\noalign{\smallskip}\hline\noalign{\smallskip}
& $\,\,d=3\,\,$ & $\, 1 \,$ &
$\, -0.542466\,$ &$\, -2.000000\,$ &
$\,-3.398381 \,$ & $\,-4.765079 \,$\\
&$\,\,d=4\,\,$ & $\,1\,$ & $\,-0.563612\,$ & $\,-2.109131\,$ & $\,-3.603718 \,$ & $\, -5.061116\,$ \\
& $\,\,d=5\,\,$ & $\,1\,$ &$\, -0.572315 \,$& $\,-2.163011\,$
& $\,-3.711951\,$ & $\,-5.216059 \,$
\\
& $\,\,d=6\,\,$ & $\,1 \,$ &$\, -0.577089 \,$& $\,-2.195673\,$
& $\, -3.780281\,$ & $\, -5.306294 \,$
\\
& $\,\,d=7\,\,$ & $\,  1 \,$  &$\, -0.580109\,$& $\,-2.217711\,$
& $\, -3.827722\,$ & $\, -5.354120  \,$
\\
& $\,\,d=8\,\,$ & $\, 1\,$ & $\,-0.582193\,$ & $\,-2.233621\,$ & $\,-3.862716\,$ & $\, -5.367078 \,$
\\[1ex]
\hline
\end{tabular}
}
\vspace{0.2cm}
\small{  \caption{ The first few eigenvalues $\lambda_n(d)$ for the (WM) case.}} \label{qnm}
\end{table}
\begin{table}[!h]
  \centering
   \scalebox{0.9}{
\begin{tabular}{|cc|ccccc|}\hline
\noalign{\smallskip}
&  & $n=0$ & $n=1$ & $n=2$ & $n=3$ & $n=4$ \\
\noalign{\smallskip}\hline\noalign{\smallskip}
& $\,\,d=3\,\,$ & $\, 1 \,$ &
$\, -0.588904 \,$ &$\, -2.181597\,$ &
$\, -3.570756 \,$ & $\, -5.043294 \,$\\
&$\,\,d=4\,\,$ & $\,1\,$ & $\,-0.623081\,$ & $\, -2.255420\,$ & $\, -3.765173 \,$ & $\, -5.326428 \,$ \\
& $\,\,d=5\,\,$ & $\,1\,$ &$\, -0.636360 \,$& $\, -2.295861\,$
& $\, -3.867878\,$ & $\, -5.479403 \,$
\\
& $\,\,d=6\,\,$ & $\,1 \,$ &$\, -0.643402 \,$& $\, -2.321548\,$
& $\, -3.932725\,$ & $\, -5.577793 \,$
\\
& $\,\,d=7\,\,$ & $\,  1 \,$  &$\, -0.647756\,$& $\, -2.339350\,$
& $\, -3.977728\,$ & $\, -5.647093  \,$
\\
& $\,\,d=8\,\,$ & $\, 1\,$ & $\, -0.650712\,$ & $\, -2.352428\,$ & $\, -4.010905\,$ & $\, -5.698800 \,$
\\[1ex]
\hline
\end{tabular}
}
\vspace{0.2cm}
\small{  \caption{ The first few eigenvalues $\lambda_n(d)$ for the (YM) case .}} \label{qnm}
\end{table}
 \section{Dynamics of blowup}
   In the lowest supercritical dimension $d=3$ Donninger proved that the linear asymptotic stability of $\phi_0$ implies its nonlinear asymptotic stability  \cite{d1,d2}. It seems feasible  that this result  can be generalized to higher dimensions $d\geq 4$ but we are not qualified to pursue this issue. We conjecture that the solution $\phi_0$ is not only nonlinearly stable but, in fact, it is a universal attractor for generic blowup. The numerical evidence for this conjecture in the case $d=3$ was given in \cite{bct,bt}; here we give an analogous evidence for $d\geq 4$.

Consider a solution $u(t,r)$ of Eq.\eqref{eq} that  starts from smooth initial data at $t=0$ and becomes singular at $(r=0,t=T)$ for some $T>0$ (recall that singularities at $r>0$ are impossible).  We want to show that generically  $u(t,r)$ converges to the self-similar solution $\phi_0(\frac{r}{T-t})$ as $t\nearrow T$
and the rate and profile of convergence are
determined by the least damped mode $e^{\lambda_1 s} v_1(y)$, that is
\begin{equation}\label{conv}
u(t,r)-\phi_0\left(\frac{r}{T-t}\right) \sim  C\, (T-t)^{-\lambda_1}\, v_1\left(\frac{r}{T-t}\right)\,,
\end{equation}
where the coefficient $C$ and blowup time $T$ depend on the initial
data.

To verify \eqref{conv}, we solved Eq.\eqref{eq} numerically for
several large, `randomly' chosen, initial data leading to blowup. In
order to keep track of the structure of the singularity developing on
progressively smaller scales, it is necessary to use an adaptive
method which refines the spatio-temporal grid near the
singularity.  Our numerical method is based on the moving mesh method
(known as MMPDE6, see \cite{hrr}) combined with the Sundman transformation, as described in
\cite{budd}, with some minor modifications and improvements specific
to the problem at hand.  This method is very efficient in
computations of self-similar singularities.

\begin{figure}[h!]
  \centering
  \includegraphics[width=0.9\linewidth]{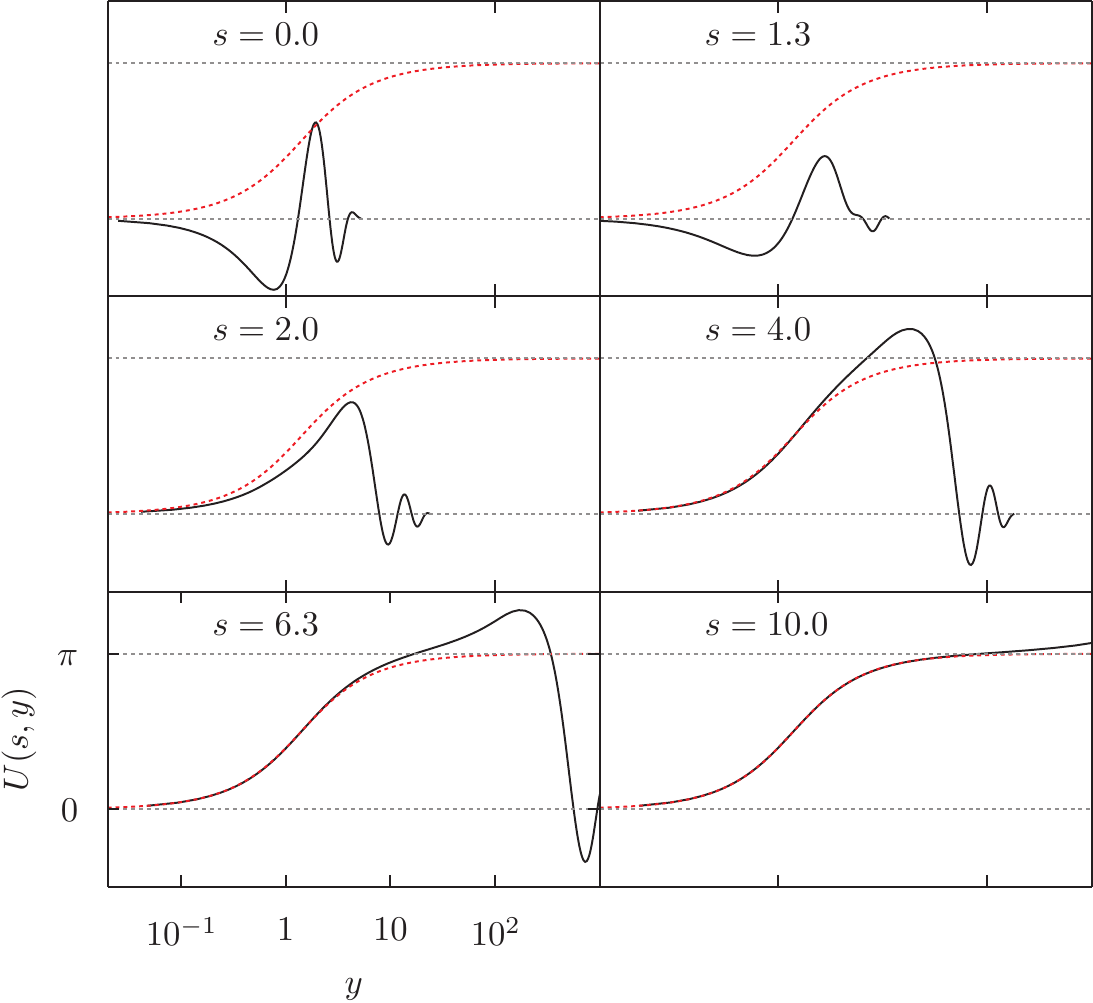}
  \caption{Snaphots from the evolution of sample initial data leading to blowup for the (WM) case in $d=4$. The solution $u(t,r)=U(s,y)$ converges to the self-similar solution
   $\phi_0(y)=2\arctan(y/\sqrt{2})$ (dashed
    red curve).}
  \label{fig:convergence}
\end{figure}

\begin{figure}[h!]
  \centering
  \includegraphics[width=\textwidth]{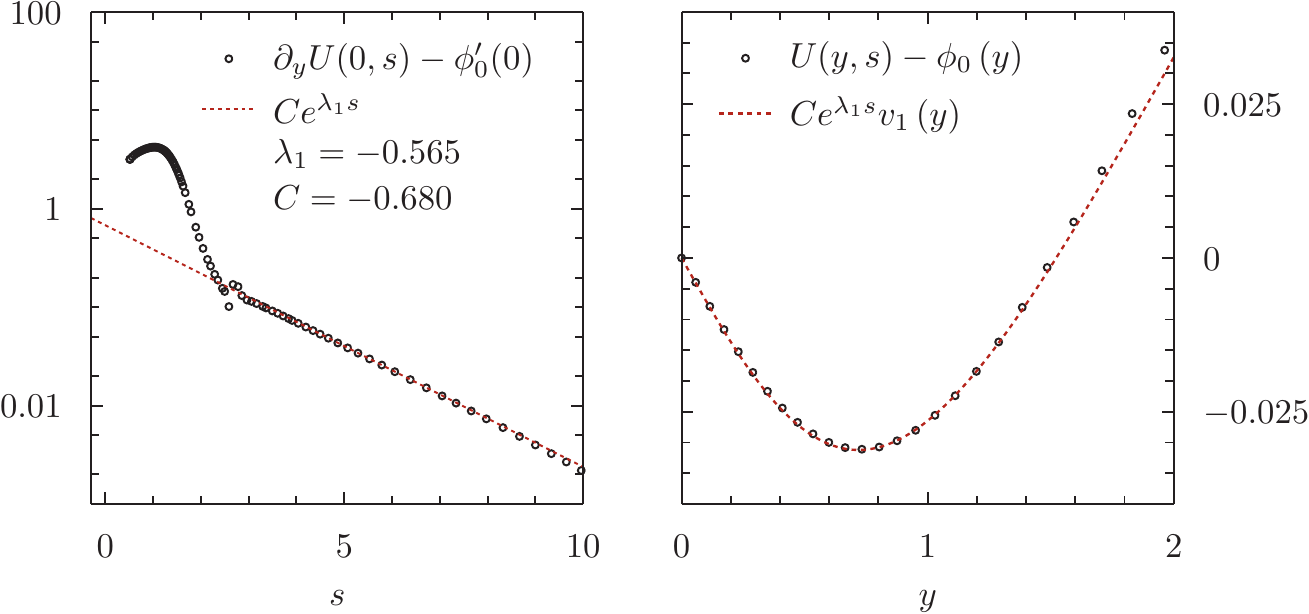}

  \caption{Numerical verification of the formula \eqref{conv} for the solution depicted in Fig.~1.  Left panel: we plot
  $(T-t) \partial_r u(t,0)-\sqrt{2}$ in  log-log scale and fit $C (T-t)^{-\lambda_1}$. The value of $\lambda_1$ obtained from the fit agrees within one percent with the value given in Table 1. Right panel: for the parameters $C$ and $\lambda_1$ determined above, the left and right hand sides of \eqref{conv} are shown to agree at time $s=4$.}
  \label{fig:rate-mode}
\end{figure}

After obtaining a numerical solution $u(r,t)$ we estimate the blow-up
time $T$ and pass to the similarity variables $s=-\log(T-t)$,
$y=r/(T-t)$ and $U(y,s)=u(r,t)$ in order to compare
analytical and numerical results.  The results are  illustrated in
Figs.~\ref{fig:convergence} and \ref{fig:rate-mode}.

\FloatBarrier
\subsubsection*{Acknowledgments.} This work was supported in part by the NCN Grant No. DEC-
2012/06/A/ST2/00397. The second author acknowledges the hospitality of the Banff International Research Station (Canada), where this work was initiated.


\begin{thebibliography}{10}

\bibitem{cst} T. Cazenave, J. Shatah, A.S. Tahvildar-Zadeh, \emph{Harmonic maps of the hyperbolic space and the development of singularities in wave maps and Yang-Mills fields,} Ann. Inst. Henri Poincar\'e 68, 315 (1998)

    \bibitem{struwe} M. Struwe, \emph{Equivariant wave maps in two space dimensions,}  Commun. Pure Appl. Math. 56(7), 815 (2003)

\bibitem{bos} P. Bizo\'n, Y.N. Ovchinnikov, I.M. Sigal, \emph{Collapse of an instanton,} Nonlinearity 17, 1179 (2004)

\bibitem{os} Y.N. Ovchinnikov, I.M. Sigal, \emph{On collapse of wave maps,} Phys. D 240, 1311 (2011)

\bibitem{rr} P. Rapha\"el, I. Rodnianski, \emph{ Stable blow up dynamics for the critical corotational wave maps and equivariant yang mills problems,} Publ. Math. Inst. Hautes Etudes Sci. 115, 1 (2012)

\bibitem{s} J. Shatah, \emph{Weak solutions and development of singularities in the $SU(2)$ $\sigma$-model,}
     Comm. Pure Appl. Math. 41, 459 (1988)

\bibitem{ts} N. Turok, D. Spergel, \emph{Global texture and the microwave background,} Phys. Rev. Lett. 64, 2736 (1990)

\bibitem{b1} P.~Bizo\'n,  \emph{Equivariant self-similar wave maps
from Minkowski spacetime into 3-sphere,} Comm. Math. Physics 215, 45 (2000)

\bibitem{b2} P. Bizo\'n, \emph{Formation of singularities in Yang-Mills equations,} Acta Phys. Polon. B 33, 1893 (2002)


\bibitem{d1} R. Donninger, \emph{On stable self-similar blowup for equivariant wave maps,} Comm. Pure Appl. Math. 64, 1095 (2011)


\bibitem{d2} R. Donninger, \emph{Stable self-similar blowup in energy supercritical Yang-Mills theory,} arXiv: 1202.1389


 \bibitem{bct} P. Bizo\'n, T. Chmaj, and Z. Tabor, \emph{Dispersion and collapse of wave maps,} Nonlinearity 13, 1411 (2000)

     \bibitem{bt} P. Bizo\'n, Z. Tabor, \emph{On blowup of Yang-Mills fields},
Phys. Rev. D 64, 121701 (2001)

\bibitem{b3} P. Bizo\'n, \emph{An unusual eigenvalue problem,} Acta. Phys. Polon. B 36, 5 (2005)

\bibitem{dsa} R. Donninger, B. Sch\"orkhuber, P.C. Aichelburg, \emph{On stable self-similar blow up for equivariant wave maps: the linearized problem,} Ann. Henri Poincar\'e 13, 103 (2012)

\bibitem{fan} H. Fan, \emph{Existence of the self-similar solutions in the heat flow of harmonic maps,} Sci. China Ser. A42,  113 (1999)

\bibitem{wein} B. Weinkove, \emph{Singularity formation in the Yang-Mills flow}, Calc. Var. Partial Differential Equations 19, 211 (2004)

\bibitem{bw} P. Bizo\'n, A. Wasserman, \emph{Non-existence of shrinkers for the harmonic map flow in higher dimensions,}  arXiv:1404.7381

\bibitem{bb} P. Biernat, P. Bizo\'n, \emph{Shrinkers, expanders, and the unique continuation
beyond generic blowup in the heat flow for harmonic
maps between spheres},     Nonlinearity 24, 2211 (2011)

   \bibitem{bier} P. Biernat, \emph{Non-self-similar blow-up in the heat flow for harmonic maps in higher dimensions,} arXiv: 1404.2209


\bibitem{nist} NIST Digital Library of Mathematical Functions,
\url{http://dlmf.nist.gov/31.2}

\bibitem{hrr} W. Huang, Y. Ren, R. D. Russell, \emph{Moving mesh partial differential equations (MMPDES) based
on the equidistribution principle,} SIAM J. Numer. Anal. 31, 709 (1994)

\bibitem{budd} C.J. Budd and J. F. Williams, \emph{How to adaptively resolve evolutionary singularities in differential equations with symmetry},
    Journal of Engineering Mathematics 66, 217-236 (2010)

\end{thebibliography}
\end{document}